\documentclass[11pt]{amsart} 
\usepackage{amsfonts,amscd}
\usepackage{psfig}
\def\N{{\mathbb N}}
\def\Z{{\mathbb Z}}

\def\R{{\mathbb R}}

\newtheorem{thr}{Theorem}[section]
\newtheorem{df}[thr]{Definition}
\newtheorem{lm}[thr]{Lemma}
\newtheorem{pr}[thr]{Proposition}

\newtheorem{cl}[thr]{Claim}

\newtheorem{thr*}{Theorem}
\newtheorem{pr*}{Proposition}
\newtheorem{co*}{Corollary}
\begin{document}
\baselineskip 22pt
\title[Nonisotopic symplectic tori in the fiber class of elliptic surfaces]{Nonisotopic symplectic
tori in the fiber class of elliptic surfaces}
\author{Stefano Vidussi}
\address{Department of Mathematics, Kansas State University, Manhattan,
KS 66506}
\email{vidussi@math.ksu.edu}
\maketitle
\baselineskip 18pt
\noindent {\bf Abstract.} The purpose of this note is to present a
construction of an
infinite family of symplectic tori $T^{p,q}$ representing an arbitrary
multiple $q[F]$ of the homology class $[F]$ of the fiber of an elliptic
surface $E(n)$, for $n \geq 3$, such that, for $i \neq j$, there is no orientation-preserving
diffeomorphism between $(E(n),T^{(i,q)})$ and $(E(n),T^{(j,q)})$. In particular, these tori are
mutually nonisotopic. This complements previous results of
Fintushel and Stern in \cite{FS2}, showing in particular the existence of
such phenomenon for a primitive class.
\section{Introduction and statement of the result} An interesting question of symplectic topology
concerns the existence, for a symplectic $4$-manifold $X$, of homologous, but not
isotopic, symplectic representatives of a given homology class. Fintushel and Stern provided, in
\cite{FS2}, the
first example of such phenomenon. Their construction, that applies to a large class of symplectic
manifolds, implies in particular that in any elliptic surface the class $2m [F]$ (where $m
\geq 2$ and
$[F]$ is the class of the elliptic fiber) can be represented by an infinite family of mutually
nonisotopic symplectic tori. Smith (\cite{S}) has been able to increase the genus of the examples,
proving that the class $2m[\Sigma_{g}]$ (where $m \geq 2$) in the (non simply-connected) surface
$\Sigma_{g}
\times S^{2}$ can be represented  by an infinite family of mutually nonisotopic symplectic curves
(whose genus can be determined by the adjunction formula).
The results above should be compared with the ones expected from a conjecture,
due to Siebert and Tian,
about the absence of such phenomena in the case of minimal rational ruled
manifolds (Siebert and Tian have in fact proven the conjecture for
several homology classes of ${\bf P}^{2}$ and $S^{2} \times S^{2}$).

These results leave open an interesting question, first pointed out
by Smith in \cite{S2}. Apart from the
problem of obtaining examples for homology classes with
odd divisibility, which appears mainly a technical question, the method
used in \cite{FS2} and \cite{S} does not allow us to obtain
nonisotopy results for primitive homology classes, as the case of
the fiber $F$ in $E(n)$.

Our purpose here is to present a different construction that produces
families of symplectic tori also in primitive homology classes, and
distinguishes their isotopy class avoiding the use of branched coverings.
This allows us to
extend (almost completely) the previous results, obtaining this way examples
of symplectic surfaces homologous but not isotopic to a complex connected
curve. Moreover, we will able to obtain a stronger result, namely that there does not
exist a orientation-preserving diffeomorphism of $E(n)$ sending one of these tori to another.
Precisely, we will prove the following \begin{thr} \label{mainth} For any $q \geq 1$ there exists an
infinite family of symplectic tori $T^{p,q}$ representing the class $q[F]$
of an elliptic surface $E(n)$, for $n \geq 3$ (where $[F]$ is the class of the fiber) such that, for
$i \neq j$, there is no orientation-preserving diffeomorphism between $(E(n),T^{(i,q)})$ and
$(E(n),T^{(j,q)})$. In particular, these tori are mutually nonisotopic. \end{thr}
We briefly sketch the argument: for each $q \geq  1$ we will consider
different homologous simple curves $K_{1}^{(p,q)}$ in the exterior of the
$3$-component link given by pushing off one component of the Hopf link.
These curves will define a family of homologous, symplectic
tori $T^{(p,q)}$ in the elliptic surface $E(n)$. We will glue copies of
the rational elliptic surface $E(1)$ along
these tori. The symplectic manifolds obtained this way are link surgery
manifolds, obtained by applying a variation of the construction of
Fintushel-Stern (introduced in \cite{V})
to a family of links introduced in Section 2. Gluing $E(1)$ along its fiber
$F$ does not depend (up to diffeomorphism of the resulting manifolds) on the
choice of the gluing map (see \cite{GS}); in particular the resulting
manifold depends only on the diffeomorphism type of the pair $(E(n),T^{(p,q)})$.
Using different tori, we will get an infinite number of mutually
nondiffeomorphic manifolds, distinguished (in a rather unusual way, see
Section \ref{many}) by the SW invariant. For two such tori $T_{1},T_{2}$
we have therefore no diffeomorphism of the pairs $(E(n),T_{1})$,
$(E(n),T_{2})$. This implies in turn that the two tori are not smoothly
isotopic.

We remark that while our examples cover cases that were
excluded in
\cite{FS2} and, {\it mutatis mutandis}, in \cite{S},
we have a price to pay, namely - as can be observed by
analyzing the construction presented in the next section - the constraint of
$n \geq 3$ of Theorem \ref{mainth} does not seem to be removable (while the
examples of \cite{FS2} exist for any elliptic surface).
\section{Construction of the family of links} 
\begin{figure}[h]
\centerline{\psfig{figure=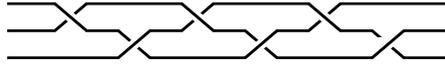,height=8mm,width=60mm,angle=0}}
\caption{ \label{borr} {\sl The braid $B_{1}$, whose closure gives the Borromean rings. }} 
\end{figure}

In this section we introduce a doubly-indexed
class of links
$\{L_{p,q},p \geq 0,q \geq 1\}$ which we will be of paramount importance in our construction:
First, denote by
$L_{0,1}$ the
$4$-component link obtained by pushing off, with respect to the $0$-framing, $2$ copies of one
component of the Hopf link,
(with the components oriented as the fibers of the Hopf fibration of $S^{3}$). Next, consider the
$3$-strand braid
$B_{1}$ of Figure \ref{borr}, and denote by $L_{1,1}$ be the $4$-component link given by the
$3$-component link
$R_{1}$ obtained by closing the braid of Figure \ref{borr}, together with the braid axis
$K_{4}$ oriented in such a way that the sublink composed by $K_{4}$ and any closed strand is the
Hopf link. 

Similarly, denote by $L_{p,1}$ the $4$-component link given by the
$3$-component link
$R_{p}$ obtained by closing the braid $B_{p}$, the composition of $p$ copies of $B_{1}$,
together with the symmetry axis $K_{4}$ oriented as before.

The link $L_{1,1}$ is the link {\it
Borromean rings plus axis}, analyzed (for different purposes) in \cite{McMT}. Its
multivariable Alexander polynomial is \begin{equation} \label{poly} \begin{array}{c}
\Delta_{L_{1,1}}(x,y,z,t) = -4 + (t + t^{-1}) + (x + x^{-1} + y + y^{-1} + z + z^{-1}) + \\ \\ -
(xy +
  (xy)^{-1} + yz + (yz)^{-1} + xz + (xz)^{-1}) + (xyz + (xyz)^{-1})
\end{array} \end{equation} where $t$ is the variable corresponding to the meridian of the axis
$K_{4}$ and $x,y,z$ correspond to the meridians of the three components given by
the closure of the strands of the braid $B_{1}$.

\begin{figure}[h]
\centerline{\psfig{figure=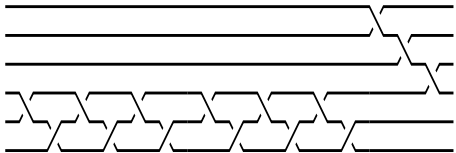,height=20mm,width=60mm,angle=0}}
\caption{ \label{itborr}
{\sl The closure of this braid with the axis gives the link $L_{2,4}$. }}
\end{figure}

The link $L_{p,q}$ is defined by modifying
$L_{p,1}$ in the following way; add, to the braid $B_{p}$, $(q-1)$ strands, which are braided to
the the first strand in the way denoted in Figure \ref{itborr}.

\begin{figure}[h]
\centerline{\psfig{figure=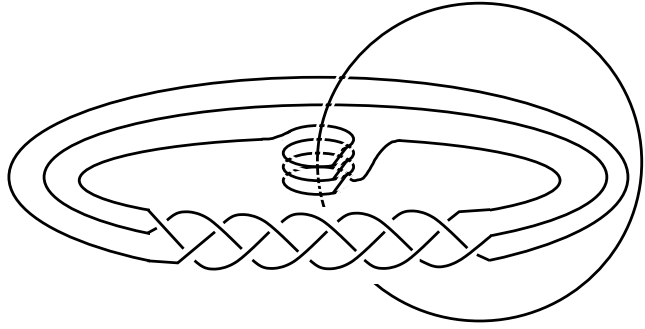,height=43mm,width=85mm,angle=0}}
\caption{ \label{stral} {\sl The link $L_{2,4}$. }}
\end{figure}

The closure of this new braid still gives
the $3$-component link $R_{p}$ (the various braids differ in fact only by Markov moves of
type II), but if we add the axis $K_{4}$, we get a new link
$L_{p,q}$, that we can visualize as obtained from $L_{p,1}$ by taking its first component
and twisting it $q$ times around $K_{4}$, see Figure \ref{stral}.

The linking matrix of $L_{p,q}$ has the form \begin{equation} \label{lima} l^{q}_{ij} =
\left( \begin{array}{cccc} - & 0 & 0 & q \\ 0 & - & 0 & 1 \\ 0 & 0 & - & 1 \\ q & 1 & 1 & -
  \end{array} \right). \end{equation}
  Observe that the linking matrix does not depend on $p$.

We will not be interested in the computation of the complete
multivariable Alexander polynomial of $L_{p,q}$; we will be content with the computation of the
reduced polynomial
$\Delta_{p,q}(s) :=
\Delta_{L_{p,q}}(s,s,s,1)$, that is determined in the following \begin{lm} \label{repo} Let
$\Delta_{p,q}(s) =
\Delta_{L_{p,q}}(s,s,s,1)$ be a specialization of the Alexander polynomial of the link $L_{p,q}$
constructed before, for $p,q \geq 1$. Then \begin{equation} \label{redpol} \Delta_{p,q}(s) =
(s^{q+2}-1)(s-1)^{3} \cdot
\prod_{j=1}^{p-1}[(1-s^{-3})(s-1)^{3} - 2(1- \cos \frac{2\pi j}{p})] \end{equation}
(with the convention that for $p=1$ the latter product is meant to be equal to $1$). \end{lm}
{\bf Proof:} To prove this equation, we need first of all the Torres formula (see e.g. \cite{T})
which in our case reads
\begin{equation} \label{torres} \Delta_{L_{p,q}}(x,y,z,1) = (x^{l^{q}_{14}}y^{l^{q}_{24}}z^{l^{q}_{34}}
-1)
\cdot \Delta_{R_{p}}(x,y,z) = (x^{q}yz-1) \cdot \Delta_{R_{p}}(x,y,z). \end{equation} where
$\Delta_{R_{p}}(x,y,z)$ is the Alexander polynomial of $R_{p}$ and the $l_{i4}^{q}$ are the
linking numbers of Equation \ref{lima}.
 To compute $\Delta_{R_{p}}(x,y,z)$,
we observe that $R_{p}$ is a periodic link, whose image under the $\Z_{p}$ action over $S^{3}$
with fixed point set the unknot $K_{4}$ is the Borromean rings $R_{1}$; from the formula for the
Alexander polynomial of periodic links (\cite{T}), and the fact that $R_{1} \cup K_{4} = L_{1,1}$,
we have
\begin{equation} \Delta_{R_{p}}(x,y,z) = \Delta_{R_{1}}(x,y,z) \cdot
\prod_{j=1}^{p-1} \Delta_{L_{1,1}}(x,y,z,\omega^{j}) \end{equation} where $\omega$ is the primitive
$p$-th root of unit. Equation \ref{poly} and explicit calculation lead then to Equation
\ref{redpol}. \qed

In the link $R_{p}$, as for the Borromean rings $R_{1}$, each component
is an unknot, and any $2$-component sublink is the trivial link. In particular, we can think at
$L_{p,q}$ as the union $K^{(p,q)}_{1}
\cup H_{3}$ where
$H_{3} = K_{2} \cup K_{3} \cup K_{4}$ is the push-off of one component of the Hopf
link (with the components $K_{2}$ and $K_{3}$ being unlinked). The links
$L_{p,q}$ - for a fixed value of $q \geq 1$ -
differ therefore from the way the unknot $K^{(p,q)}_{1}$ is linked to the $3$-component
link
$H_{3}$. In particular, if we consider the link exterior $S^{3} \setminus \nu H_{3}$, the link
exteriors $S^{3} \setminus \nu L_{p,q}$ are obtained by removing nonisotopic circles. The case
$L_{0,1}$ corresponds to the removal of the circle $K_{1}^{(0,1)}$ isotopic to $\mu(K_{4})$, the
meridian of $K_{4}$. In the 
case of $L_{p,q}$ instead we are removing the circle $K_{1}^{(p,q)}$
which is homologous to $q \mu(K_{4})$ in $H_{1}(S^{3} \setminus \nu H_{3})$, as
from the linking matrix of Equation \ref{lima} we deduce that $K_{1}^{(p,q)}$ has linking number $q$ 
with the axis $K_{4}$, and $0$ with the other two components. In what follows we will consider the
circle $K_{1}^{(p,q)}$, as well as any other link component,
endowed of the framing defined by a spanning disk.
\section{Link surgery manifolds associated to $L_{p,q}$} \label{tori}
In this section we will construct the family of $4$-manifolds used to prove Theorem \ref{mainth}.
We start by recalling briefly the definition of link
surgery manifold (see \cite{FS}), in the modified form introduced in \cite{V}. Consider an
$m$-component link $K
\subset S^{3}$ and take an homology basis of simple curves $(\alpha_{i},\beta_{i})$ of intersection
$1$ in the boundary of the link exterior. Next, take $m$ elliptic surfaces $E(n_{i})$ and define
the manifold
\begin{equation} \label{lisu} E(K) = (\coprod E(n_{i}) \setminus \nu F_{i}) \cup_{F_{i}
\times S^{1} = S^{1} \times
\alpha_{i} \times \beta_{i}} (S^{1} \times (S^{3} \setminus \nu K)), \end{equation} where the
orientation reversing diffeomorphism between the boundary $3$-tori is defined so to identify
$F_{i}$ with $S^{1} \times \alpha_{i}$ and acts as complex conjugation  on the remaining circle factor.

It is well known that in general the fiber sum above is not
well defined and, for a fixed choice of homology basis, the smooth structure of the manifold above
could depend on various choices, but because of the use of elliptic surfaces the manifold we will
discuss will not be affected by this indeterminacy.

We have now a simple claim, whose
proof follows by the definition of the elliptic surface
$E(n)$ as an iterated fiber sum of elliptic surfaces. Fix $\{n_{i}\} = \{1,1,n-2\}$.
\begin{cl}
\label{coel} Let
$H_{3}$ be the $3$-component link obtained by pushing off one copy of one component of the
Hopf link; then we can chose the homology basis $(\alpha_{i},\beta_{i})$ so that $E(H_{3}) = E(n)$.
\end{cl} {\bf Proof}:
This claim follows from the observation that
\begin{equation} S^{1} \times (S^{3} \setminus \nu H_{3}) = T^{2} \times (S^{2} \setminus \nu
\{p_{2},p_{3},p_{4}\}) \end{equation} so that choosing  $(\alpha_{i},\beta_{i}) =
(\lambda(K_{i}),-\mu(K_{i}))$ for $i=2,3$ and $(\alpha_{4},\beta_{4}) =
(\mu(K_{4}),\lambda(K_{4}))$ we have an explicit presentation of $E(n)$.
\qed

In $E(n)$
defined as above, the image of the class of the curve $\mu(K_{4})$ under the injective map
\begin{equation} \label{injmap} H_{1}(S^{3} \setminus \nu H_{3},\Z) \stackrel{S^{1} \times
(\cdot)}{\longrightarrow} H_{2}(S^{1} \times    (S^{3} \setminus \nu H_{3}),\Z) \longrightarrow
H_{2}(E(n),\Z) \end{equation} is the class of the elliptic fiber. More precisely, the image of
the torus $S^{1} \times \mu(K_{4})$ in $E(n)$ is identified with a copy of the elliptic fiber $F$.

Consider now the images $T^{(p,q)}$ of the tori $S^{1}
\times K_{1}^{(p,q)}$ under the injection \begin{equation} S^{1} \times (S^{3} \setminus \nu
H_{3}) \hookrightarrow E(n);
\end{equation} these compose a family of embedded, self-intersection zero framed tori. 
We have the following \begin{pr} Up to isotopy, the tori $T^{(p,q)}$ are
symplectic submanifolds of $E(n)$, homologous to $qF$, where $F$ is the fiber of the elliptic
fibration.
\end{pr} {\bf Proof:}  The statement on homology follows from the fact that the circles
$K_{1}^{(p,q)}$ are all homologous to $q\mu(K_{4})$ in $H_{1}(S^{3} \setminus \nu H_{3},\Z)$, and
the class
$[T^{(p,q)}]$ coincides therefore
with the image of $q[\mu(K_{4})]$ under the map of Equation \ref{injmap},
i.e. it is the multiple  $q[F]$ of the class of the fiber.

In order to prove that the $T^{(p,q)}$ are
symplectic, we will present $E(n)$, together
with its symplectic structure, as a symplectic fiber sum in the following way:
we perform a surgery with coefficients respectively $\infty,\infty,0$ to
$K_{2} \cup K_{3} \cup K_{4} \subset S^{3}$ (i.e. ultimately a $0$-surgery to the unknot $K_{4}
\subset S^{3}$) to obtain the three manifold $S^{1} \times S^{2}$, in which
the cores $C_{i}$ of the solid tori used in the surgery (specifically $K_{2}$ and $K_{3}$ itself,
plus a curve isotopic to $\mu(K_{4})$) are framed, essential curves, whose
framing induces one for the tori $S^{1} \times C_{i} \subset S^{1} \times S^{1} \times S^{2}$.
Then we have \begin{equation} E(n) = \coprod_{i=2}^{4} E(n_{i}) \#_{F_{i} = S^{1}\times
C_{i}} S^{1} \times (S^{1} \times S^{2}). \end{equation}
Note that, by the definition of fiber sum and because of the framings of $S^{1} \times C_{i}$,
this construction coincides with the one of Claim \ref{coel}.

In $S^{1} \times S^{2}$ the curves
$C_{i}$ are transverse to the fiber $S^{2}$ of the obvious fibration (which extends the
$D^{2}$ fibration of $S^{3} \setminus
\nu K_{4} = S^{1} \times D^{2}$) and if we denote by $\phi \in \Omega^{1} (S^{1} \times S^{2},\R)$ a closed
nondegenerate representative of that fibration, for any curve $C$ in $S^{3} \setminus \nu K_{4}$
which is transversal to the disk fibration, we have pointwise $\phi(C) > 0$; as a consequence,
endowing $S^{1} \times S^{1} \times S^{2}$ of the symplectic structure $\phi \wedge dt + \epsilon
\psi$ (with $\psi$ a volume form on the fiber $S^{2}$ and $\epsilon$ sufficiently small),
the tori $S^{1} \times
C_{i}$ (more generally, any torus $S^{1} \times C$ as above) are symplectic in $S^{1} \times
S^{1} \times S^{2}$ and, consequently, in $E(n)$. The curves $K_{1}^{(p,q)} \subset S^{3}
\setminus
\nu K_{4}$ are (up to isotopy) transverse to the disk fibration, and the tori $T^{(p,q)}$ are
therefore symplectic. \qed

We can now introduce the link surgery manifolds associated to the links $L_{p,q}$.
These are defined as in Equation \ref{lisu}, but we can also present them as fiber sum of $E(n)$ and
$E(1)$ along the embedded tori $T^{(p,q)} \subset E(n)$ and $F \subset E(1)$. This is the content
of  the next definition, in which we write also the Seiberg-Witten polynomial of the manifold. Fix
$(n_{1},n_{2},n_{3},n_{4}) = (1,1,1,n-2)$:
\begin{df} \label{fourma} Let $L_{p,q}$ be the
$4$-component link  considered above, and define \begin{equation} \label{asma}
E(L_{p,q}) = (\coprod_{i=1}^{4} E(n_{i}) \setminus \nu F_{i}) \cup_{F_{i} \times S^{1} = S^{1}
\times
\alpha_{i} \times \beta_{i}} (S^{1} \times (S^{3} \setminus \nu
L_{p,q})) = E(n) \#_{T^{(p,q)} = F} E(1) \end{equation} where $(\alpha_{1},\beta_{1}) =
(\lambda(K^{(p,q)}_{1}),-\mu(K^{(p,q)}_{1}))$,
$(\alpha_{i},\beta_{i}) =
(\lambda(K_{i}),-\mu(K_{i}))$ for $i = 2,3$ and $(\alpha_{4},\beta_{4}) =
(\mu(K_{4}),\lambda(K_{4}))$. The SW polynomial is given by the product of the relative SW
invariants
\begin{equation}
\label{swpoly} SW(E(L_{p,q})) = (\prod_{i=1}^{4} SW(E(n_{i}) \setminus \nu F_{i})) \cdot SW(S^{1} \times (S^{3} \setminus \nu
L_{p,q})) =
(t-t^{-1})^{n-3} \Delta^{s}_{L_{p,q}}(x^{2},y^{2},z^{2},t^{2}) \end{equation} where
$\Delta^{s}$ is the symmetrized version of the multivariable Alexander polynomial.
\end{df} The latter statement follows from Theorem 2.7 of \cite{Ta} (see also \cite{FS}),
as the homology class of the fiber of
$E(n-2)$ (the elliptic surface glued to $S^{1} \times $(axis of $L_{p,q}$)) is identified with the
image of
$S^{1} \times \mu(K_{4})$ in $E(L_{p,q})$.

Note that, although we made explicit a choice of
curves $(\alpha_{1},\beta_{1})$ in Definition \ref{fourma}, the smooth structure of the resulting 
manifold is independent of such choice, i.e. depends ultimately only on the
diffeomorphism type of $(E(n),T^{(p,q)})$.

For sake of notation, we will omit reference to the number $n$ for the manifold in Equation
\ref{asma}, its value being clear from the context. We observe that $E(L_{0,1})$ is just
$E(n+1)$ (see Claim \ref{coel}), while $E(L_{1,1})$, for $n=3$, is the interesting manifold
considered in \cite{McMT}, with the  presentation discussed in \cite{V}.
\section{Infinitely many nonisotopic tori} \label{many} In this section we will prove Theorem
\ref{mainth}, namely we will show that, for a fixed value of
$q$, there are infinitely many diffeomorphism types of pairs $(E(n),T^{(p,q)})$.
In order to prove that two tori $T^{(i,q)}$, $T^{(j,q)} \subset E(n)$ define different pairs for
$i \neq j$ it would be sufficient to prove that the manifolds $E(L_{i,q})$, $E(L_{j,q})$ have
different SW polynomial. This means that there does not exist any automorphism of the manifold,
inducing an automorphism of the second cohomology group  which sends $SW(E(L_{i,q}))$ to
$SW(E(L_{j,q}))$ (note that, when comparing the SW polynomials of two manifolds, as the ones
appearing in Equation
\ref{swpoly}, we must consider the fact that the variables with the same symbol could refer to
different cohomology classes for the two manifolds). Proving such a result appears to be
quite a challenging problem (also considering the fact that we do not have a complete
knowledge of the SW polynomials of our manifolds).

We will not attempt here to prove this, and we will limit ourselves to the proof of a weaker
statement, that is anyhow sufficient to prove the statement of Theorem
\ref{mainth}. The model of proof we will exploit here could find application also in other similar
problems, where the explicit comparison of SW polynomials is difficult.

We will start, for sake of example, to work out in detail (and with a proof
which differs from the general case) the case  of two preferred
tori, among the ones defined in Section \ref{tori}, namely $T^{(0,q)}$ and
$T^{(1,q)}$. The proof that these tori define different pairs
 constitutes, in some sense, a ``finite'' version of Theorem \ref{mainth}.
To obtain such a result, we will use in a rather weak way SW theory,
building from the following observation: Let $d(X)$ be the dimension of the
the vector subspace of $H^{2}(X,\R)$ spanned by SW basic classes of $X$;
then $d(X)$ is a smooth invariant of $X$. We use this fact to prove the
following \begin{thr} For any $q \geq 1$ the manifolds $E(L_{0,q})$ and
$E(L_{1,q})$ are nondiffeomorphic (symplectic) manifolds. \end{thr}
{\bf Proof:} in order to prove that, we will show that
$d(E(L_{0,q})) = 2$ while $d(E(L_{1,q})) > 2$.
The first statement follows from the explicit
computation of the Alexander polynomial of $L_{0,q}$: we can observe that
$L_{0,q}$ is a graph link obtained by connected sum along $K_{4}^{*}$
of a $2$-component link given by the unknot $K^{I}_{4}$ and its $(1,q)$
cable $K_{1}^{(0,q)}$ with the $3$-component link given by the unknot
$K^{II}_{4}$ and two copies $K_{2} \cup K_{3}$ of the meridian. We leave to
the reader the application of the results of \cite{EN} to verify
that the Alexander polynomial of this graph link is \begin{equation}
\Delta_{L_{0,q}}(x,y,z,t) = (t-1)^{2} \frac{x^{q}t^{q} - 1}{xt-1}.
\end{equation} (For similar computations see e.g. \cite{V2}.)
In particular, this polynomial depends on only two variables,
and the nonzero terms span a $2$-dimensional subspace of
$H_{1}(S^{3} \setminus \nu L_{0,1},\R)$. From this and Equation \ref{swpoly}
the statement about $d(E(L_{0,q}))$ follows.
For $d(E(L_{1,q}))$, we can observe that the span of the
nonzero terms of $\Delta_{L_{1,q}}(x,y,z,t)$ is bounded by
below by the span of nonzero terms of the reduced polynomial
$\Delta_{L_{1,q}}(x,y,z,1)$ which is
given, according to Equation \ref{torres}, by
\begin{equation} \Delta_{L_{1,q}}(x,y,z,1) = (x^{q}yz - 1)(x-1)(y-1)(z-1).
\end{equation} The span of nonzero terms of this polynomial,
as is easily verified, has dimension $3$; using Equation \ref{swpoly} again
we obtain that $d(E(L_{1,q})) \geq 3$ (note that the fact that the SW
polynomial reduced at $t=1$ is zero, for $n > 3$, does not affect this).
This completes the proof. \qed

We will discuss now the general case.
We will prove the following \begin{thr} For any $q \geq 1$ the family
$\{E(L_{p,q})\}_{p \in \N}$
contains an infinite number of nondiffeomorphic (symplectic) manifolds. \end{thr} {\bf Proof:} To
prove this statement it is sufficient to prove that, if we denote by $\beta_{p}$ the number of basic
classes of the manifold $E(L_{p,q})$ (for a fixed $q$), we have $\lim_{p} \beta_{p} = + \infty$.
We will start by proving this for the case of $n=3$, where the $SW$ invariant ``coincides" with
the Alexander polynomial of $L_{p,q}$, as written in Equation \ref{swpoly}.
In this case we can observe that the number of basic classes of $E(L_{p,q})$ coincides with the
number of nonzero terms in
$\Delta_{L_{p,q}}(x,y,z,t)$. Such a number is bounded by below by the number $\tau_{p}$ of
nonzero terms in the reduced polynomial $\Delta_{p,q}(s)$ of Lemma \ref{repo}, that we rewrite
here by convenience: \begin{equation} \label{rewri} \Delta_{p,q}(s) = \sum_{k} a_{p,k}s^{k} =
(s^{q+2}-1)(s-1)^{3} \cdot
\prod_{j=1}^{p-1}[(1-s^{-3})(s-1)^{3} - 2(1- \cos \frac{2\pi j}{p})]. \end{equation}
In order to estimate $\tau_{p}$ we observe that the number of nonzero terms $a_{p,k}$ of a Laurent
polynomial in $s$ satisfies the inequality of Lemma
\ref{maxle} in the appendix, i.e. $\tau_{p} \geq \frac{1}{2} \rho_{p} + 1$ where $\rho_{p}$ is the
number of nonzero real roots of $\Delta_{p,q}$. The proof that $\lim_{p} \rho_{p} = + \infty$ will
therefore prove our statement. It follows from elementary arguments
that the equation $(1 - s^{-3})(s-1)^{3} = 2(1-\cos \alpha)$ has exactly $2$ real reciprocal solutions
$0 < s_{1}(\alpha) < 1 < s_{2}(\alpha)$ for $0 < \alpha \leq \pi$, which differ for different values of
$\alpha$. As a consequence each of the first $[\frac{p-1}{2}]$ factors appearing in the product of
Equation \ref{rewri} contributes two roots to
$\rho_{p}$, and we have \begin{equation} \label{fineq} \rho_{p} \geq  1 + 2 [ \frac{p-1}{2} ].
\end{equation}  This proves the statement for $n=3$.

We point out that the estimate
on the number of terms is not optimal; in particular for odd $q$ it is not difficult to prove that
$\tau_{p} = 6p+1$.

To prove the statement for $n > 3$ we consider the
specialization of the SW polynomial given by \begin{equation} SW_{p}(s,s,s,t) = (t-t^{-1})^{n-3}
\Delta^{s}_{L_{p,q}}(s^{2},s^{2},s^{2},t^{2}) \end{equation}
Once again, to prove that $\lim_{p} \beta_{p} = + \infty$ it is sufficient to prove that the number of nonzero
terms in $SW_{p}(s,s,s,t)$ goes to infinity with $p$. We can rewrite such a
two-variable polynomial as \begin{equation} SW_{p}(s,s,s,t) =: \sum_{k} (t-t^{-1})^{n-3}
a_{p,k}(t)s^{k}
\end{equation} where, in the last identity, we define $a_{p,k}(t)$ as the polynomial in $t$
that appears in writing $\Delta^{s}_{L_{p,q}}(s^{2},s^{2},s^{2},t^{2})$ as a power series in $s$.
If we consider the number ${\tilde \tau_{p}}$ of nonzero coefficients $(t-t^{-1})^{n-3}a_{p,k}(t)$,
this is bounded by below by the number of nonzero $a_{p,k}(1)$; but the set of the latter coefficients
(with a reparametrization for $k$ that takes account of the symmetrization and the ``squaring'' of the $s$ variable) coincides the set of the coefficients $a_{p,k}$ of Equation
\ref{rewri}: therefore ${\tilde \tau}_{p} \geq \tau_{p}$ and Equation
\ref{fineq} asserts that this number diverges with
$p$. This completes the proof of the statement. \qed

Notice that, although $\beta_{i} \neq \beta_{j}$ implies $E(L_{i,q}) \neq
E(L_{j,q})$, the condition ${\tilde \tau_{i}} \neq {\tilde \tau_{j}}$
is instead not sufficient to prove this, as we
cannot guarantee that the specializations of the Alexander polynomials are the same.

As the family of manifolds obtained by gluing $E(1)$ to $E(n)$ along different $T^{(p,q)}$, for a
fixed $q$, contains infinitely many nondiffeomorphic manifolds, infinitely many pairs
$(E(n),T^{(p,q)})$ are not diffeomorphic. In particular there are infinitely many
nonisotopic symplectic tori $T^{(p,q)}$. This completes the proof of Theorem \ref{mainth}.
\section{Appendix}
In this Appendix we give a proof of the Lemma used in Section \ref{many}. (It is likely
that this statement already exists in literature, but we have not been able to find a reference). We thank Maximilian Seifert for suggesting us the proof of this Lemma.
\begin{lm} \label{maxle} Let $p(z)$ be a nontrivial real Laurent polynomial. Denote by $\rho$ the
number of nonzero real roots (counted without multiplicity) and by $\tau$ the number of terms of the
polynomial. Then we have the inequality $\rho \leq 2\tau -2$. \end{lm} {\bf Proof:} Assume first
that $p(z)$ is an ordinary polynomial of degree $n$ satisfying \begin{equation} \label{condi} p(z)
= \sum_{k=0}^{n} a_{k}z^{k}, \ \ a_{n} \neq 0, \ \  a_{0} \neq 0 \end{equation} and denote by
$\gamma$ the number of holes appearing in the polynomial plus
$1$, where we define by {\it hole} a string of consecutive powers $z^{d},z^{d-1},...,z^{d-*}$ with
coefficient equal to zero and $1 \leq d < n$ (e.g. $p(z) = 2z^{6}-4z^{2}+3$ has $\gamma = 3$). By
obvious reasons, $\gamma \leq \tau$. Introduce now the family of integer pairs
$(n_{l},m_{l})_{l=1,...,\gamma}$, with $n_{l} \geq m_{l} > n_{l+1} + 1$, defined in such a way that
\begin{equation} p(z) =
\sum_{l=1}^{\gamma}
\sum_{k=m_{l}}^{n_{l}} a_{k} z^{k}; \end{equation} this means that $a_{d} \neq 0 \leftrightarrow
d \in [m_{l},n_{l}]$ for some $1 \leq l \leq \gamma$.

We will first prove, by
induction over $\gamma$, that for a polynomial as in Equation \ref{condi} we have the inequality 
\begin{equation} \rho \leq
\sum_{l=1}^{\gamma} (n_{l} - m_{l}) + 2\gamma -2. \end{equation} This inequality is trivially true for
$\gamma = 1$. Assume by inductive hypothesis that it holds true for $\gamma-1$: we want to
prove it for
$\gamma$. Take the first $(n_{\gamma} + 1)$ derivatives of $p(z)$ and denote \begin{equation} 
\label{gine}
\begin{array}{cc}
q(z) := (\frac{d}{dz})^{n_{\gamma}+1} p(z)  = \sum_{l=1}^{\gamma-1} \sum_{k=m_{l}}^{n_{l}}
a_{k}
\frac{k!}{(k-(n_{\gamma}+1))!}z^{k-(n_{\gamma}+1)} = 
\\ \\  = z^{m_{\gamma-1} - (n_{\gamma}+1)} \sum_{l=1}^{\gamma -1}
\sum_{k=m_{l}}^{n_{l}} a_{k}
\frac{k!}{(k-(n_{\gamma}+1))!}z^{k-m_{\gamma-1}}  =:  z^{m_{\gamma-1} - (n_{\gamma}+1) }
{\tilde q(z)}. \end{array} \end{equation} The polynomial ${\tilde q(z)}$ has one hole less
than $p(z)$ and satisfies the conditions of Equation \ref{condi}: we can thus apply the
inductive hypothesis for it. Moreover, the roots of
$q(z)$ coincide with the roots of ${\tilde q(z)}$, plus the root $z=0$: in particular we have
\begin{equation} \label{inez} \rho(q(z)) = \rho({\tilde q(z)}) \leq \sum_{l=1}^{\gamma-1} (n_{l}
- m_{l}) +2 \gamma - 4. 
\end{equation} By Rolle's theorem, the number of real zeroes of $p(z)$ is bounded in terms of
 the zeroes of
its derivative: more precisely we have, from Equation \ref{inez} and the fact that $m_{\gamma} =
0$, the inequality
\begin{equation} \rho(p(z)) \leq
\rho(q(z)) + 1 + (n_{\gamma} + 1) \leq \sum_{l=1}^{\gamma-1} (n_{l} - m_{l}) + n_{\gamma} + 2\gamma
- 2 = \sum_{l=1}^{\gamma} (n_{l} - m_{l}) + 2 \gamma - 2
\end{equation} which is what we wanted to prove. This completes our induction.

Now we can observe that $\tau =
\sum_{l=1}^{\gamma}(n_{l} - (m_{l} - 1))$. Applying this to Equation \ref{gine}, together 
with the inequality $\gamma \leq \tau$,
proves the Lemma when $p(z)$ is an ordinary polynomial. The statement for a general Laurent
polynomials is readily obtained from this, by multiplying the polynomial with a suitable power of
$z$ in order to get an ordinary polynomial of the form of Equation \ref{condi}. \qed
\vspace*{5mm}

{\bf Acknowledgements:} I would like to thank Maximilian Seifert for 
the proof of Lemma \ref{maxle}. I would like to thank also 
Patrick Brosnan and Ki Heon Yun for discussions.

\end{document}